\documentclass[a4paper,12pt,leqno]{article}

\usepackage{amssymb,amsmath}
\usepackage{xypic}
\xyoption{all}
\usepackage{graphicx}

\newtheorem{defn}{Definition}[section]

\newtheorem{corollary}[defn]{Corollary}
\newtheorem{rem}[defn]{Remark}
\newtheorem{exm}[defn]{Example}
\newtheorem{lemma}[defn]{Lemma}
\newtheorem{theorem}[defn]{Theorem}
\newtheorem{notat}[defn]{Notation}
\newtheorem{newpar}[defn]{}

\newtheorem{xdefn}{Definition.}
\newtheorem{xproposition}{Proposition.}
\newtheorem{xcorollary}{Corollary.}
\newtheorem{xrem}{Remark.}
\newtheorem{xexm}{Example.}
\newtheorem{xlemma}{Lemma.}
\newtheorem{xtheorem}{Theorem.}
\newtheorem{xnotat}{Notation.}
\newtheorem{xnewpar}{\it}
\newtheorem{xproof}{{\it Proof. }}
\newtheorem{xproofof}{{\it Proof}}

\newenvironment{definition}{\begin{defn}\em}{\end{defn}}

\newenvironment{example}{\begin{exm}\em}{\end{exm}}

\newenvironment{proof}{\begin{xproof}\em}{\end{xproof}}

\newenvironment{newparagraph*}[1]{\begin{xnewpar}\hspace*{-1.5mm}{#1}. \rm}{\end{xnewpar}}

\newenvironment{definition*}{\begin{xdefn}\em}{\end{xdefn}}
\newenvironment{remark*}{\begin{xrem}\em}{\end{xrem}}
\newenvironment{example*}{\begin{xexm}\em}{\end{xexm}}
\newenvironment{notation*}{\begin{xnotat}\em}{\end{xnotat}}
\newenvironment{proposition*}{\begin{xproposition}}{\end{xproposition}}
\newenvironment{corollary*}{\begin{xcorollary}}{\end{xcorollary}}
\newenvironment{lemma*}{\begin{xlemma}}{\end{xlemma}}
\newenvironment{theorem*}{\begin{xtheorem}}{\end{xtheorem}}

\newcommand\RR{\mathbb{R}}
\def\qed{\hspace{0.3cm}{\rule{1ex}{2ex}}}
\newcommand\V{\bigvee}
\newcommand\ie{i.e.}
\newcommand\eg{e.g.}
\newcommand\st{\mid}
\newcommand\cf{\textrm{cf.}}
\newcommand\downsegment{{\downarrow}}

\newcommand\spp{\varsigma}

\newcommand\topology{\operatorname{\Omega}}

\newcommand\Mod{\mathbf{Mod}}
\newcommand\Frm{\mathbf{Frm}}
\newcommand\Loc{\mathbf{Loc}}
\newcommand\Etale{\mathbf{LH}}

\newcommand\opp[1]{{#1}^{\textrm{op}}}

\newcommand\ident{\mathrm{id}}

\newcommand\inner[2]{{\langle #1,#2 \rangle}}

\newcommand\hide[1]{}

\newcommand\sets{\mathbf{Sets}}
\newcommand\sections{\mathit{\Gamma}}

\newcommand\hmb{\mathbf{HMB}}
\newcommand\sh{\mathbf{Sh}}

\newcommand\Mat{\mathbf{Mat}}

\newenvironment{eq}{\setcounter{equation}{\arabic{defn}}\begin{equation}}{\end{equation}\setcounter{defn}{\arabic{equation}}}
\newenvironment{eqarray}{\setcounter{equation}{\arabic{defn}}\begin{eqnarray}}{\end{eqnarray}\setcounter{defn}{\arabic{equation}}}

\begin{document}

\title{Sheaves as modules\thanks{Research supported in part by the Funda\c{c}\~ao para a Ci\^encia e a Tecnologia through the Program POCI 2010/FEDER, namely via Centro de An\'{a}lise Matem\'{a}tica, Geometria e Sistemas Din\^{a}micos (first author), Centro de \'{A}lgebra da Universidade de Lisboa (second author), and the grant POCI/MAT/55958/2004.}}

\author{\sc Pedro Resende and Elias Rodrigues}

\date{}

\maketitle

\begin{abstract}
We revisit sheaves on locales by placing them in the context of the theory of quantale modules. The local homeomorphisms $p:X\to B$ are identified with the Hilbert $B$-modules that are equipped with a natural notion of basis. The homomorphisms of these modules are necessarily adjointable, and the resulting self-dual category yields a description of the equivalence between local homeomorphisms and sheaves whereby morphisms of sheaves arise as the ``operator adjoints'' of the maps of local homeomorphisms.
\end{abstract}


\section{Introduction}\label{introduction}

Quantales, at least those with certain properties, can be regarded as point-free generalized (``non-commutative'') spaces \cite{Mulvey,PR,Rosenthal}. In particular, they subsume localic \'etale groupoids \cite{AIM} and, to some extent, C*-algebras \cite{KR}.
There are several proposals for what could be an appropriate notion of sheaf on such a space (see, \eg, \cite{BB,G00,G01,MN,JoelPhD,Rosenthal2,IsarPhD,S07,W}), and the present paper has grown out of an effort to understand the relations between them and in particular to find concrete examples.
Such sheaves can be identified with suitable quantale modules in a way that at least models equivariant sheaves on \'etale groupoids \cite{am}, but we believe that even just for sheaves on locales the theory is interesting and worth being presented separately. This is the purpose of the present paper, which moreover provides useful background for \cite{am}.

\subsection*{\sc Overview}

Sheaves on a locale $B$ (or on a topological space) can be described in several different ways. Besides the two classical formulations, namely as local homeomorphisms $p:X\to B$ and as separated and complete presheaves $F:\opp B\to\sets$, they can also be presented as ``$B$-valued sets''. These are (possibly infinite) projection matrices $ E:\sections\times\sections\to B$ with entries taken from $B$, where $\sections$ can be identified with the set of local sections of a sheaf in the usual sense. This idea, with the locale $B$ playing the role of a commutative ring, goes back to \cite{FS} (expositions of it can be found in \cite[\S2]{Borceux} and \cite[pp.\ 502--513]{Elephant}), and a variant is obtained from \cite{W} by taking the ``ring'' to be the quantaloid that arises as the Cauchy completion of $B$ (see also \cite[pp.\ 77--82]{Rosenthal2}). Approaching sheaves in the language of quantale modules provides us with another way, which in addition is quite natural and in some sense helps unify the others.

In this paper we shall begin, in \S\ref{sec:mapsasmods}, with a description in terms of $B$-modules of the open locale maps, and even arbitrary maps, into $B$, and we shall obtain direct translations between local homeomorphisms and locale modules, early on identifying two isomorphic ``categories of sheaves as modules'', $B$-$\Etale$ and $B$-$\sh$, with the same objects. The morphisms of $B$-$\Etale$, called \emph{maps}, coincide with the maps between local homeomorphisms, whereas the morphisms of $B$-$\sh$ are called \emph{sheaf homomorphisms} because they correspond directly to natural transformations between sheaves.

Local homeomorphisms, or sheaves, turn out to correspond to particular instances of the Hilbert $Q$-modules of \cite{Paseka}, which are an analogue of Hilbert C*-modules \cite{Lance} where C*-algebras are replaced by involutive quantales.
Hence, a Hilbert $Q$-module is a $Q$-module equipped with a $Q$-valued ``inner product''. In this paper $Q$ will be a locale $B$, and, as we shall see in \S\ref{sec:hilbmods}, such a module corresponds to a sheaf precisely when it has a subset playing, in the quantale context, the role of a Hilbert basis --- the basis ``vectors'' are local sections of the sheaf.
The construction of a $B$-valued set from such a $B$-module is immediate: it is the inner product restricted to basis elements (the ``metric'' of the inner product).

We shall see in \S\ref{sec:adjmaps} that for Hilbert $B$-modules equipped with Hilbert bases the $B$-module homomorphisms are precisely the functions $h$ that have adjoints $h^\dagger$ defined in terms of the inner product in the same way as one does for Hilbert spaces:
\[\langle h(x),y\rangle = \langle x,h^\dagger(y)\rangle\;.\]
Hence, we have a category $B$-$\hmb$ of such modules which is equipped with a strong self-duality $(-)^\dagger:\opp{(B\textrm{-}\hmb)}\to B\textrm{-}\hmb$ that, in particular, yields the isomorphism between $B$-$\Etale$ and $B$-$\sh$ as a restriction. The adjoint $(f^*)^\dagger$ coincides with the left adjoint $f_!$ of $f^*$, thus providing us with an example of a situation where categorical adjoints coincide with ``operator adjoints''.

\subsection*{\sc Notation and  terminology}

We shall use fairly standard notation and terminology. In particular, the word \emph{locale} is used as a synonym for \emph{frame};
a \emph{homomorphism} of locales $h:X\to Y$ is a function that preserves arbitrary joins (including the least element $0$) and finite meets (including the greatest element $1$);
and a \emph{continuous map} of locales (or simply a \emph{map}) $f:Y\to X$ is a homomorphism in the opposite direction, $f^*:X\to Y$, also referred to as the \emph{inverse image homomorphism} of $f$.

The category of locales and their maps is denoted by $\Loc$ and it is referred to as the \emph{category of locales}.

If $B$ is a locale then by a \emph{$B$-module} is meant, as usual, a complete lattice $X$ equipped with a sup-preserving action $B\otimes X\to X$ of the commutative monoid $(B,\wedge,1)$ --- we use direct sum and tensor product notation for locales and modules in analogy with the notation for commutative rings and their modules, as in \cite{JT}. 

Throughout the paper, $B$ is a fixed but arbitrary locale.

\subsection*{\sc Acknowlegment}

The authors thank the referee for comments to an earlier version of this paper.

\section{Continuous maps as modules}\label{sec:mapsasmods}\label{sec:shhoms1}

\subsection*{\sc General continuous maps}

Let $p:X\to B$ be a map of locales. Then $X$ is a $B$-module by ``change of base ring'' along the homomorphism $p^*:B\to X$:  the action is given by, for all $x\in X$ and $b\in B$,
\[b x=p^*(b)\wedge x\;.\]
It follows that $b1=p^*(b)$ and thus this module satisfies the condition
\begin{eq}\label{Bloc}
b x=b1\wedge x\;,
\end{eq}%
which, as we shall see, completely characterizes the modules that arise in this way.
(This condition has been called \emph{stability} in \cite{RV}, in the more general context of modules over unital quantales.)
We remark that the action of such a module distributes over meets of non-empty sets $S\subset X$ in the right variable:
$b\left(\bigwedge S\right) = b1\wedge\bigwedge S
=\bigwedge_{x\in S} (b1\wedge x) = \bigwedge_{x\in S} bx$.

Let us define some terminology:

\begin{definition}\label{def:Blocales}
Let $B$ be a locale. By a \emph{$B$-locale} will be meant a locale $X$ equipped with a structure of $B$-module satisfying (\ref{Bloc}). A \emph{homomorphism} of $B$-locales is a homomorphism of locales that is also a homomorphism of $B$-modules, and a \emph{map} $f:X\to Y$ of $B$-locales is defined to be a homomorphism $f^*:Y\to X$ of $B$-locales. The \emph{category of $B$-locales}, denoted by $B$-$\Loc$, has as objects the $B$-locales and as morphisms the maps of $B$-locales. We shall denote the category $\opp{(B\textrm{-}\Loc)}$ by $B$-$\Frm$ (the \emph{category of $B$-frames}).
\end{definition}

\begin{theorem}\label{thm:isomorphism}
The category $B$-$\Loc$ is isomorphic to $\Loc/B$.
\end{theorem}

\begin{proof}
Each object $p:X\to B$ of $\Loc/B$ gives us a $B$-locale, as we have seen in the beginning of this section. Conversely, let $X$ be a $B$-locale. Define a function
$\phi:B\to X$
by
\[\phi(b)=b1\;.\]
We have $\phi(1)=11=1$, $\phi(b\wedge c)=(b\wedge c)1 = b(c1)=b1\wedge c1=\phi(b)\wedge \phi(c)$,
and $\phi(\V_\alpha b_\alpha)=(\V_\alpha b_\alpha)1=\V_\alpha b_\alpha 1=\V_\alpha \phi(b_\alpha)$; that is, $\phi$ is a homomorphism of locales, and thus we have obtained a map $p:X\to B$ defined by $p^*=\phi$.
This correspondence between objects of $\Loc/B$ and $B$-locales is clearly a bijection.

In order to see that the categories are isomorphic let $p:X\to B$ and $q:Y\to B$ be objects of $\Loc/B$, and let
$f:X\to Y$ be a map of locales. We show that $f$ is a morphism from $p$ to $q$ in $\Loc/B$ if and only if it is a map from $X$ to $Y$ in $B$-$\Loc$. Let $b\in B$ and $y\in Y$. We have
\[f^*(by)=f^*(q^*(b)\wedge y)=(q\circ f)^*\wedge f^*(y)\]
and also
\[bf^*(y)=p^*(b)\wedge f^*(y)\;.\]
It follows that if $p=q\circ f$ then $f^*(by)=bf^*(y)$ for all $b\in B$ and $y\in Y$; that is, if $f$ is in $\Loc/B$ then $f^*$ is a homomorphism of $B$-modules and thus $f$ is in $B$-$\Loc$. Conversely, if $f^*$ is a homomorphism of $B$-modules then letting $y=1$ above we obtain $q\circ f=p$. \qed
\end{proof}

From now on we shall freely identify $B$-modules with their associated locale maps, for instance calling $B$-module to a map $p:X\to B$, and for convenience we shall often refer to $p$ as the \emph{projection} of the $B$-module.

\subsection*{\sc Open maps}

\begin{definition}
A $B$-locale $X$ is \emph{open} if its projection $p$ is an open map of locales; that is, $p^*$ has a left adjoint $p_!$ which is a homomorphism of $B$-modules (but not in general a homomorphism of $B$-locales).
\end{definition}

It is obvious that the direct image $p_!$ of the projection $p$ of an open $B$-locale satisfies the property
\[p_!(x)x=x\;,\]
for $p_!(x)x=p^*(p_!(x))\wedge x$ and thus the equality $p_!(x)x=x$ is equivalent to the unit of the adjunction $p_!\dashv p^*$.
This has a converse: if $\spp:X\to B$ is $B$-equivariant and monotone and it satisfies
\[\spp(x)1\ge x\]
then $\spp$ is left adjoint to the map $(-)1:B\to X$; the condition $\spp(x)1\ge x$ is the unit of the adjunction and the counit $\spp(b1)\le b$ is an immediate consequence of the equivariance, for $\spp(b1)=b\land\spp (1)\le b$. This actually holds for any $B$-module, but for a $B$-locale $X$ the condition $\spp(x)1\ge x$ (equivalently, $\spp(x)x=x$ because $\spp(x)x=\spp(x)1\land x$) implies that $X$ is open. Summarizing, we have:

\begin{theorem}\label{thm:openmaps}
A $B$-locale $X$ is open if and only if there is a monotone equivariant map
\[\spp:X\to B\]
such that the following (necessarily equivalent) conditions are satisfied for all $x\in X$:
\begin{eqarray}
\spp(x)1 &\ge& x\\
\spp(x)x &\ge& x\\
\spp(x)x &=& x\;.\label{strongsupp}
\end{eqarray}%
Furthermore there is at most one such map $\spp$. If one exists it is necessarily a $B$-module homomorphism and it coincides with the direct image $p_!$ of the projection $p$ of $X$.
\end{theorem}

Alternatively, a slightly different characterization of open $B$-locales is the following:

\begin{theorem}
Let $X$ be a locale which is also a $B$-module (but not necessarily a $B$-locale). Then $X$ is an open $B$-locale if and only if there is a monotone equivariant map
$\spp:X\to B$
such that (\ref{strongsupp}) holds for all $x\in X$.
\end{theorem}

\begin{proof}
By the previous theorem any open $B$-locale has such a map $\spp$, so we only have to prove the converse.
Assume that $\spp:X\to B$ is monotone, equivariant, and that it satisfies (\ref{strongsupp}). Then we have
\[b1\land x = \spp(b1\land x)(b1\land x)\le\spp(b1)x=(b\land\spp(1))x\le bx\;.\]
The converse inequality, $bx\le b1\land x$, is obvious, and thus (\ref{Bloc}) holds. This shows that $X$ is an open $B$-locale. \qed
\end{proof}

If $X$ is an open $B$-locale with projection $p$ then $p_!$ will be referred to as the \emph{support} of $X$, and we shall usually write $\spp$ instead of $p_!$, following the analogous notation for supported quantales \cite{AIM}. Similarly, we may refer to $p_!(x)$ as the \emph{support} of $x$.

\begin{example}\label{example:freemodule}
Let $S$ be any set. Then the free $B$-module generated by $S$, which is the function module $B^S$ of maps $f:S\to B$, is an open $B$-locale whose support is defined by $\spp(f)=\V_{s\in S} f(s)$. The projection of the $B$-locale is the obvious map $p:\coprod_{s\in S} B\to B$, where $\coprod_{s\in S} B$ is the coproduct in $\Loc$ of as many copies of $B$ as there are elements in $S$; in other words, $p^*:B\to B^S$ is the diagonal homomorphism that to each $b\in B$ assigns the map $f:S\to B$ such that $f(s)=b$ for all $s\in S$.
\end{example}

\begin{example}\label{exm:omegalocales}
In a topos, any locale $X$ has a unique $\Omega$-locale structure determined by the continuous map $!_X:X\to\Omega$, and $X$ is an open $\Omega$-locale precisely if it is open in the usual sense \cite{JT}. The tensor product $B\otimes X$ (the product $B\times X$ in $\Loc$) is a $B$-locale with action $a(b\otimes x)=(a\wedge b)\otimes x$ and projection $\pi_1^*(b)=b\otimes 1$, and it is open if $X$ is open (because $\pi_1$ is the pullback of $!_X$ along $!_B$). Its support is computed from the $\Omega$-action on $B$ by $\spp(b\otimes x)=\spp(x)b$. In $\sets$ every locale is open and the support of $B\otimes X$ is defined by the conditions $\spp(b\otimes 0)=0$ and $\spp(b\otimes x)=b$ if $x\neq 0$.
\end{example}

\subsection*{\sc Local homeomorphisms}

Let $p:X\to B$ be a local homeomorphism, and let $\sections$ be a cover of $X$ (\ie, $\sections\subset X$ and $\V\sections=1$) such that, on each open sublocale determined by an element of $\sections$, $p$ restricts to a homeomorphism onto its image; that is, for each $s\in\sections$ there is a commutative square
\begin{eq}\label{fact}
\vcenter{\xymatrix{
B\ar[rr]^{p^*}\ar@{->>}[d]_{(-)\land \spp(s)} & & X\ar@{->>}[d]^{(-)\land s}\\
\downsegment \spp(s)\ar[rr]^{\theta_s} & & \downsegment s }}
\end{eq}%
such that $\theta_s$ is an isomorphism. Then we have, for each $b\le \spp(s)$,
\[\theta_s(b)=p^*(b)\land s = bs\]
and
\[\spp(\theta_s(b))=\spp(bs)=b\land\spp(s)=b\;.\]
Hence, the restriction to $\downsegment s$ of $\spp$ splits $\theta_s$, and thus it coincides with $\theta_s^{-1}$. This motivates the following definition:

\begin{definition}\label{def:localsection}
Let $X$ be an open $B$-locale. A \emph{local section} of $X$ is an element $s\in X$ such that for all $x\le s$ we have
\begin{eq}
\spp(x)s=x\label{iso2}\;.
\end{eq}%
The set of local sections of $X$ is denoted by $\sections_X$, and $X$ is defined to be \emph{\'etale} if
$\V\sections_X=1$.
\end{definition}

If $s$ is a local section then
\[p^*(1)\land s=1\land s=s=p^*(\spp(s))\land s\;,\]
which means that the homomorphism $((-)\land s)\circ p^*$ factors as in (\ref{fact}). Moreover, the equivariance of $\spp$ gives us $\spp(bs)=b$ for all $b\le\spp(s)$; this, together with (\ref{iso2}), ensures that $\theta_s$ an isomorphism, and therefore we have the following characterization of local homeomorphisms:

\begin{theorem}
An open $B$-locale is \'{e}tale if and only if its projection is a local homeomorphism.
\end{theorem}

Let us denote by $\Etale$ the subcategory of $\Loc$ whose objects are the locales and whose morphisms are the local homeomorphisms. It follows from a basic property of local homeomorphisms that $\Etale/B$ is a full subcategory of $\Loc/B$, and thus $\Etale/B$ is isomorphic to the following category, which provides our first example of a ``category of sheaves as modules'':

\begin{definition}
The \emph{category of \'etale $B$-locales}, denoted by $B$-$\Etale$, is the full subcategory of $B$-$\Loc$ whose objects are the \'etale $B$-locales.
\end{definition}

\subsection*{\sc Sheaf homomorphisms}

Although $B$-$\Etale$ is meant to be a ``category of sheaves as modules'', it is not a category of modules; that is, $\opp{(B\textrm{-}\Etale)}$, rather than $B$-$\Etale$, is a subcategory of $B$-$\Mod$. In order to remedy this let us first introduce the following terminology:

\begin{definition}\label{def:shhom}
Let $X$ and $Y$ be \'etale $B$-locales. By a \emph{sheaf homomorphism}
\[h:X\to Y\]
will be meant a homomorphism of $B$-modules satisfying the following two conditions:
\begin{enumerate}
\item $h(\sections_X)\subset \sections_Y$;
\item $\spp(h(s))=\spp(s)$ for all $s\in\sections_X$ (equivalently, $\spp(h(x))=\spp(x)$ for all $x\in X$).
\end{enumerate}
The sheaf homomorphisms form a category that we shall denote by $B$-$\sh$.
\end{definition}

The motivation for this terminology comes from the fact that, denoting by $\sh(B)$ the category of sheaves on $B$ in the usual sense (sheaves are separated and complete presheaves $\opp B\to\sets$ and their morphisms are the natural transformations), we have a functor 
$G:B\textrm{-}\sh\to\sh(B)$
(which is part of an equivalence of categories --- see the comments below) such
that: (i) $G$ assigns to each \'etale $B$-locale $X$ the sheaf $G_X:\opp B\to\sets$ defined for each $b\in B$ by
$G_X(b)=\{s\in\sections_X\st \spp(s)=b\}$,
with the restriction map $G_X(b)\to G_X(a)$ for each pair $a\le b$ in $B$ being given by $s\mapsto as$; (ii) $G$ assigns to each sheaf homomorphism $h:X\to Y$ the natural transformation $G_h:G_X\to G_Y$ whose components $(G_h)_b$ are all defined by $s\mapsto h(s)$.

\begin{lemma}
If $f:X\to Y$ is a map of \'etale $B$-locales then $f_!:X\to Y$ is a sheaf homomorphism.
\end{lemma}

\begin{proof}
Let $f:X\to Y$ be a map of $B$-locales with projections $p$ and $q$, respectively. Then $f_!$ satisfies $\spp_Y(f_!(x))=q_!(f_!(x))=p_!(x)=\spp_X(x)$ for all $x\in X$. In addition, $f_!$ is $Y$-equivariant, and thus it is $B$-equivariant for the module structures of $X$ and $Y$ induced by $f^*\circ q^*$ and $q^*$, respectively. Finally, composing $f$ with a local section of $p$ yields a local section of $q$ --- a module theoretic proof of this is as follows: if $y\le f_!(s)$ and $s\in\sections_X$ then
\[y=y\wedge f_!(s) = f_!(f^*(y)\wedge s) = f_!(\spp(f^*(y)\wedge s)s)
=\spp(f^*(y)\wedge s)f_!(s)\;,\]
and thus $y$ is a restriction of $f_!(s)$. \qed
\end{proof}

Hence, we have a functor $\mathcal S:B\textrm{-}\Etale\to B\textrm{-}\sh$ which is the identity on objects and to each map $f$ assigns $f_!$. Using the (localic) correspondence between sheaves and local homeomorphisms (as in, \eg, \cite[\S2]{Borceux} or \cite[pp.\ 502--513]{Elephant}) it is not hard to see that $\mathcal S$
is part of an adjoint equivalence of categories whose other functor is the composition
\[B\textrm{-}\sh\stackrel G\longrightarrow \sh(B)\stackrel\Lambda\longrightarrow \Etale/B\stackrel\cong\longrightarrow B\textrm{-}\Etale\;,\]
where, concretely, $\Lambda$ can be the functor that to each sheaf assigns its locale of closed subobjects as in \cite[\S2.2]{Borceux}. But in fact one can prove something stronger:

\begin{theorem}\label{thm:postponed}
The functor $\mathcal S:B\textrm{-}\Etale\to B\textrm{-}\sh$ is an isomorphism.
\end{theorem}

A direct proof of this, using properties of Hilbert modules, will be postponed until \S\ref{sec:shhoms}.

\section{Hilbert modules}\label{sec:hilbmods}

\subsection*{\sc Basic definitions and properties}

Hilbert $Q$-modules are an analogue of Hilbert C*-modules where C*-algebras are replaced by quantales. They have been studied by Paseka mainly as a means of importing results and techniques from operator theory into the context of quantales (see, \eg, \cite{Paseka2}), and also in connection with theoretical computer science \cite{Paseka}.
We begin by recalling
this notion
in the special case that interests us in this paper, namely when the involutive quantale $Q$ is the locale $B$.

\begin{definition}
By a \emph{pre-Hilbert $B$-module} will be meant a $B$-module $X$ equipped with a function
\[\inner --:X\times X\to B\]
called the \emph{inner product}, that satisfies the following axioms, for all $x,y\in X$ and $b\in B$:
\begin{eqarray}
\inner{bx}y&=&b\land\inner x y\\
\left\langle\V_\alpha x_\alpha,y\right\rangle &=&\V_\alpha\inner{x_\alpha}y\\
\inner x y&=&\inner y x\;.
\end{eqarray}%
(In short, a symmetric $B$-valued ``bilinear'' form.)
A \emph{Hilbert $B$-module} is a pre-Hilbert $B$-module whose inner product is non-degenerate,
\begin{eq}
\inner x -=\inner y -\Rightarrow x=y\;,
\end{eq}%
and it is said to be \emph{strict} (``positive definite'') if it satisfies
\[\inner x x=0\Rightarrow x=0\;.\]
\end{definition}

A useful consequence of non-de\-ge\-ne\-racy is the following:

\begin{lemma}\label{lem:stability}
Let $X$ be a Hilbert $B$-module. Then for all $b\in B$ and $x\in X$ we have
\[bx = b1\wedge x\;.\]
Hence, in particular, if $X$ is a locale it is a $B$-locale.
\end{lemma}

\begin{proof}
The inequality $bx\le b1\wedge x$ is immediate. For the other, it suffices to show that for all $y\in X$ we have $\inner{b1\wedge x}y\le\inner{bx}y$:
\[\inner{b1\wedge x}y\le\inner{b1}y\wedge\inner x y=b\wedge\inner 1 y\wedge \inner x y=
\inner{bx}y\;. \qed\]
\end{proof}

\subsection*{\sc Supported modules and open $B$-locales}

Now we see some relations between Hilbert $B$-modules and open $B$-locales.

\begin{definition}
A (pre-)Hilbert $B$-module $X$ is \emph{supported} if it satisfies the condition
\[\inner x x x=x\]
for all $x\in X$. (Hence, $X$ is strict.)
\end{definition}

\begin{theorem}\label{thm:SuppHilbIsOpen}
Any locale $X$ which is also a supported Hilbert $B$-module is an open $B$-locale; its support function $\spp$ is defined by $\spp x=\langle x,x\rangle$.
\end{theorem}

\begin{proof}
Assume that $X$ is a locale equipped with a structure of Hilbert $B$-module. Then it is a $B$-locale due to \ref{lem:stability}. Besides, the function $\spp:X\to B$ defined by
\[\spp x = \langle x,x\rangle\]
is monotone and $B$-equivariant, and by hypothesis it satisfies $\spp (x) x = x$, whence by \ref{thm:openmaps} $X$ is open. \qed
\end{proof}

There is a partial converse to this theorem:

\begin{theorem}\label{thm:weaklynondeg}
Let $X$ be an open $B$-locale. Then $X$ is a supported pre-Hilbert $B$-module whose inner product is \emph{weakly non-degenerate} in the sense that if
$\langle x,z\rangle = \langle y,z\rangle$
for all $z\in X$ then $\neg x=\neg y$ (where $\neg$ is the Heyting algebra pseudo-complement: $\neg x=x\rightarrow 0$).
\end{theorem}

\begin{proof}
If $X$ is an open $B$-locale we define
$\langle x,y\rangle = \spp(x\wedge y)$.
Being a $B$-locale implies that for all $x,y\in X$ and all $b\in B$ we have
\[bx\wedge y = (b1\wedge x)\wedge y = b1\wedge(x\wedge y) = b(x\wedge y)\;,\]
and thus
\[\langle bx,y\rangle = \spp(bx\wedge y) = \spp(b(x\wedge y)) = b\wedge\spp(x\wedge y) = b\wedge\langle x,y\rangle\;.\]
If $y\in X$ and $(x_\alpha)$ is a family of elements in $X$ we have
\[\left\langle\V_\alpha x_\alpha,y\right\rangle = \spp\left(\V_\alpha x_\alpha\wedge y\right) = \V_\alpha\spp(x_\alpha\wedge y) = \V_\alpha\langle x_\alpha,y\rangle\;.\]
Since $\langle-,-\rangle$ is of course symmetric, it follows that $X$ is a pre-Hilbert $B$-locale.
For the weak non-degeneracy let $x,y\in X$ be such that $\spp(x\wedge z)=\spp(y\wedge z)$ for all $z\in X$. Then, letting $z=\neg y$, we obtain
\[0=\spp 0 = \spp(y\wedge \neg y)=\spp(x\wedge \neg y)\;,\]
and thus $x\wedge\neg y=\spp(x\wedge\neg y)(x\wedge\neg y)=0$. Hence, $\neg y\le\neg x$. Similarly, letting $z=\neg x$ we conclude that $\neg x\le\neg y$. \qed
\end{proof}

The ``weakly'' in the theorem cannot be dropped. In order to see this, consider as an example of open map the first projection $\pi_1:\RR^2\to\RR$. Let $U$ be an open ball centered on $(0,0)\in\RR^2$, and let $V=U\setminus\{(0,0)\}$. For all open sets $W\in\topology(\RR^2)$ we have $\pi_1(U\cap W)=\pi_1(V\cap W)$, but $U\neq V$ and thus the inner product associated to $\pi_1$ is degenerate.

\subsection*{\sc Hilbert bases}

Let us introduce a natural notion in the context of Hilbert $B$-modules, namely the analogue of a Hilbert basis of a Hilbert space. As we shall see, the existence of such a basis has strong consequences, notably modules equipped with a Hilbert basis are necessarily \'etale $B$-locales.

\begin{definition}
Let $X$ be a pre-Hilbert $B$-module. By a \emph{Hilbert basis} of $X$ is meant a subset $\sections\subset X$ such that for all $x\in X$ we have
\[x = \V_{s\in \sections} \inner x s s\;.\]
(In particular, $\sections$ is therefore a set of $B$-module generators for $X$.)
\end{definition}

A Hilbert basis in this sense is not an actual basis as in linear algebra because there is no freeness (we only have projectivity --- see \ref{usefullemma} below). Therefore one might be better off calling it a Hilbert system of generators, but for the sake of simplicity we shall retain the shorter terminology.

\begin{example}\label{exm:freeHilbertbasis}
Let $S$ be a set. The free $B$-module $B^S$ (\cf\ Example \ref{example:freemodule}) has a Hilbert basis $\sections$ consisting of the ``unit vectors'' $f^{(s)}:S\to B$; for each $s\in S$ we define $f^{(s)}=\iota_s(1)$ where $\iota_s:B\to B^S\cong\bigoplus_{s\in S} B$ is the coproduct injection corresponding to the $s$-labeled copy of $B$. This definition of $f^{(s)}$ makes sense in any topos and it is equivalent, if $S$ is decidable, to the following:
\[f^{(s)}(t) = \left\{\begin{array}{ll}
1 & \textrm{if }t = s\\
0 & \textrm{if }t \neq s \;.
\end{array}\right.\]
\end{example}

The existence of a Hilbert basis has many useful consequences. In particular, any pre-Hilbert $B$-module with a Hilbert basis is necessarily supported, hence strict, and it is projective:

\begin{lemma}\label{usefullemma}
Let $X$ be a pre-Hilbert $B$-module and let $\sections\subset X$. If $\sections$ is a Hilbert basis then the following properties hold, for all $x,y\in X$.
\begin{enumerate}
\item\label{ulproj} $X$ is a projective $B$-module.
\item\label{ul1} $\V\sections = 1$. ($\sections$ is a \emph{cover} of $X$.)
\item\label{ul2} If $\inner x s=\inner y s$ for all $s\in \sections$ then $x=y$. (Hence, $X$ is a Hilbert module.)
\item\label{ul3} $\inner x y=\V_{s\in\sections}\inner x s\land\inner s y$.
\item\label{ul4}\label{usefullemma-5} $\inner x x x=x$. (Hence, $X$ is supported.)
\item\label{ul5}\label{bilin4} $\inner x y\le\inner x x$.
\item\label{ul6}\label{usefullemma-4} For all $s\in\sections$ the following conditions are equivalent:
\begin{enumerate}
\item\label{ula} $x\le s$;
\item\label{ulb} $x=\inner x x s$;
\item\label{ulc} $x=\inner x s s$.
\end{enumerate}
\item\label{ulmatrix} The $B$-valued matrix $M:\sections\times\sections\to B$ defined by $m_{st}=\inner s t$ (the ``metric'' of the  inner product) is a projection matrix: $M^T=M^2=M$ (hence, $M$ defines a $B$-set in the sense of \cite{FS} --- in fact one can see that it is a complete $B$-set).
\end{enumerate}
Conversely, $\sections$ is a Hilbert basis if $\inner--$ is non-degenerate and \ref{ul3} holds.
\end{lemma}

\begin{proof}
Assume that $\sections$ is a Hilbert basis. The first eight properties are proved as follows.
\begin{enumerate}
\item[{\rm \ref{ulproj}.}]
Since $\sections$ is a set of $B$-module generators and $B^{\sections}$ is a free module, there is a quotient of $B$-modules $\varphi:B^\sections\to X$ given by $\varphi(f)=\V_{s\in\sections} f(s)s$, and in the opposite direction we define another homomorphism of $B$-modules $\psi:X\to B^{\sections}$ by $\psi(x)(s) = \inner x s$. This splits $\varphi$, showing that $X$ is a retract of a free module:
\[\varphi(\psi(x))=\V_{s\in\sections} \psi(x)(s)s=\V_{s\in\sections} \inner x s s= x\;.\]
\item[{\rm \ref{ul1}.}]
$1=\V_{s\in\sections}\inner 1 s s\le\V_{s\in\sections}1 s=\V\sections$.
\item[{\rm \ref{ul2}.}] If $\inner x s=\inner y s$ for all $s\in \sections$ then
$x=\V_{s\in\sections}\inner x s s=\V_{s\in\sections}\inner y s s=y$.
\item[{\rm \ref{ul3}.}]
$\inner x y=\left\langle{\V_{s\in\sections}\inner x s s},y\right\rangle=\V_{s\in\sections}\inner x s\land\inner s y$.
\item[{\rm \ref{ul4}.}] For all $x\in X$ and $s\in\sections$ we have
$\inner{\inner x x x}s=\inner x x\land\inner x s
=\V_{t\in\sections}\inner x t\land\inner t x\land\inner x s
=\V_{t\in\sections}\inner x t\land\inner x t\land\inner x s=\inner x s$,
and thus by the non-degeneracy we conclude $\inner x x x=x$.
\item[{\rm \ref{ul5}.}] Using \ref{ul4} we have $\inner x y=\inner{\inner x x x}y=\inner x x\land\inner x y$.
\item[{\rm \ref{ul6}.}] Either of the equations \ref{ulb} or \ref{ulc} implies \ref{ula}, of course, so let us assume that $x\le s$ in order to verify the converse implication. By \ref{ul5} we have
$\inner x x=\inner x s$ and thus \ref{ulb} and \ref{ulc} are equivalent; in addition, we have $\inner x x s\ge\inner x x x=x$, and, conversely,
$\inner x x s=\inner x s s\le\V_{t\in\sections}\inner x t t=x$, whence $x=\inner x x s=\inner x s s$.
\item[{\rm \ref{ulmatrix}.}]
We have $M=M^T$ by definition of the inner product, and $M=M^2$ follows from \ref{ul3}.
\end{enumerate}
Now assume that $\inner - -$ is non-degenerate and that \ref{ul3} holds. Then for all $x,y\in X$ we have
\[\left\langle{\V_{s\in\sections}\inner x s s},y\right\rangle=\V_{s\in\sections}\inner x s\land\inner s y=\inner x y\;,\]
and by the non-degeneracy we obtain $\V_{s\in\sections}\inner x s s=x$. \qed
\end{proof}

Proposition \ref{usefullemma}-\ref{ulmatrix} has a converse, namely every projection matrix has an associated Hilbert module (in fact a locale --- \cf\ \ref{lem:HMBisLOC}) with a Hilbert basis (we shall write $Mf$ for the product of the matrix $M$ by the ``column vector'' $f:S\to B$ --- that is, writing also $f_s$ instead of $f(s)$ for such ``vectors'', we have $(Mf)_s = \V_{t\in S} m_{st}\land f_t$):

\begin{lemma}\label{lemma:matrixtomod}
Let $S$ be a set and $M:S\times S\to B$ a $B$-valued projection matrix. Then the subset of $B^S$
\[M B^S=\{Mf\st f\in B^S\}\]
is a Hilbert $B$-module with the same inner product as $B^S$, it has a Hilbert basis $\sections$ consisting of the functions $\tilde s:S\to B$ defined, for each $s\in S$, by $\tilde s_t=m_{ts}$ ($\tilde s$ is the ``$s^\textrm{th}$-column'' of $M$), and  for all $s,t\in S$ we have
\[m_{st}=\inner{\tilde s}{\tilde t}\;.\]
\end{lemma}

\begin{proof}
The assignment $j:f\mapsto Mf$ is a $B$-module endomorphism of $B^S$, and $MB^S$ is its image, hence a submodule of $B^S$. Next note that $\sections$ is a subset of $MB^S$ because for each $t\in S$ we have
$\tilde t=M\tilde t\in MB^S$:
\[\tilde t_s=m_{st}=(M^2)_{st}=\V_{u\in S} m_{su}\land m_{ut} = \V_{u\in S} m_{su}\land\tilde t_u =(M\tilde t)_s \;.\]
For all $f\in MB^S$ we have $f=Mf$ and thus it follows that, for all $s\in\sections$,
\[\inner f{\tilde s}=\V_t f_t\land\tilde s_t=\V_t f_t\land m_{ts}=
\V_t m_{st}\land f_t=(Mf)_s=f_s\;.\]
Hence, $\sections$ is a Hilbert basis because for all $s\in \sections$ we have
\[\left(\V_t\inner f{\tilde t}\tilde t\right)_s = \left(\V_t f_t\tilde t\right)_s
=\V_t f_t\land \tilde t_s = \V_t m_{st}\land f_t = (Mf)_s=f_s\;. \qed\]
\end{proof}

\subsection*{\sc \'Etale $B$-locales}

Now we establish an equivalence between local homeomorphisms, on one hand, and Hilbert $B$-modules equipped with Hilbert bases, on the other.

\begin{lemma}\label{lem:HMBisLOC}
Any Hilbert $B$-module with a Hilbert basis is necessarily a $B$-locale and it arises, up to isomorphism, as in \ref{lemma:matrixtomod}.
\end{lemma}

\begin{proof}
Let $X$ be a Hilbert $B$-module with a Hilbert basis $\sections$, let $M$ be the matrix determined by $m_{st}=\inner s t$ for all $s,t\in\sections$, and let $\varphi:B^\sections\to X$ be the $B$-module quotient defined by $\varphi(f) = \V_{s\in \sections} f_s s$. Recalling that the inner product is non-degenerate we have,  for all $f,g\in B^\sections$, the following series of equivalences:
\begin{eqnarray*}
\varphi(f)=\varphi(g) &\iff& \forall_{t\in\sections}\ \inner{\varphi(f)}t=\inner{\varphi(g)}t\\
&\iff&\forall_{t\in\sections}\ \left\langle\V_{s\in \sections} f_s s,t\right\rangle = \left\langle\V_{s\in \sections} g_s s,t\right\rangle\\
&\iff& \forall_{t\in\sections}\ \V_{s\in \sections} f_s\land\inner s t=\V_{s\in \sections} g_s\land\inner s t\\
&\iff& \forall_{t\in\sections}\  (Mf)_t=(Mg)_t\\
&\iff& Mf=Mg\;.
\end{eqnarray*}
This shows that the $B$-module surjection $\varphi$ factors uniquely through the quotient $f\mapsto Mf : B^\sections\to MB^\sections$ and an isomorphism of $B$-modules $X\stackrel\cong\to MB^\sections$. 
Finally, in order to conclude that $X$ is a $B$-locale it suffices to show that it is a locale, due to \ref{lem:stability}, or, equivalently, that $MB^\sections$ is a locale. Consider the $B$-module endomorphism $j:f\mapsto Mf$ of $B^\sections$, as in the proof of \ref{lemma:matrixtomod}. It is easy to prove directly that $MB^\sections$ is a locale (the restriction of $j$ to $\downsegment(M 1_{B^\sections})$ is a closure operator whose fixed points define a subframe --- not a sublocale --- of $\downsegment(M 1_{B^\sections})$), but in fact this is already known, for $MB^\sections$ coincides with the set of $B$-subsets of $M$ as in \cite[Def.\ 2.8.9 and Prop.\ 2.8.11]{Borceux}. \qed
\end{proof}

\begin{theorem}
Let $X$ be a $B$-module. The following conditions are equivalent:
\begin{enumerate}
\item $X$ can be equipped with a structure of pre-Hilbert $B$-module for which there is a Hilbert basis;
\item $X$ can be equipped with a structure of Hilbert $B$-module for which there is a Hilbert basis;
\item $X$ is an \'etale $B$-locale.
\end{enumerate}
\end{theorem}

\begin{proof}
The first two conditions are equivalent due to \ref{usefullemma}-\ref{ul2}.
Let us prove that 2 implies 3. Let $X$ be a Hilbert $B$-module (and thus also a $B$-locale, by \ref{lem:HMBisLOC}). By \ref{usefullemma}-\ref{ul4} $X$ is supported, and thus by \ref{thm:SuppHilbIsOpen} it is an open $B$-locale with $\spp$ defined by $\spp x=\inner x x$.
Now let $s\in\sections$, and let $x\le s$. We have $\spp(x)s=x$, by \ref{usefullemma}-\ref{ul6}, and thus $s$ is a local section in the sense of \ref{def:localsection}. Hence, $\sections\subset\sections_X$. Since by \ref{usefullemma}-\ref{ul1} we know that $\V\sections=1$, we conclude that $X$ is \'{e}tale.

Now let us prove that 3 implies 1. If $X$ is \'{e}tale it is open and thus by \ref{thm:weaklynondeg} it is a supported pre-Hilbert $B$-module with the inner product defined by
$\inner x y=\spp(x\land y)$.
For each $x\in X$ we have
\[x=1\land x=\left(\V_{s\in\sections_X} s\right)\land x
=\V_{s\in\sections_X} s\land x\;,\]
and, by the definition of local section, $s\land x=\spp(s\land x)s=\inner x s s$, thus showing that $\sections_X$ is a Hilbert basis. \qed
\end{proof}


This provides us with an analogy between the view of sheaves on $B$ as Hilbert $B$-modules, on one hand, and vector bundles on a compact space $X$ as Hilbert $C(X)$-modules, on the other. The analogy extends to the fact that, just as vector bundles on $X$ are projective $C(X)$-modules, sheaves on $B$ are projective $B$-modules (due to \ref{usefullemma}-\ref{ulproj}). But, of course, the analogue of Swan's theorem does not hold because projective $B$-locales (\ie, the $B$-locales that are projective objects in $B$-$\Mod$) are not necessarily \'etale $B$-locales. For instance, if $B=\Omega$ this is just the statement that not every locale which is projective as a sup-lattice is necessarily a free sup-lattice. The following diagram, where $f$ is a retraction, illustrates this:
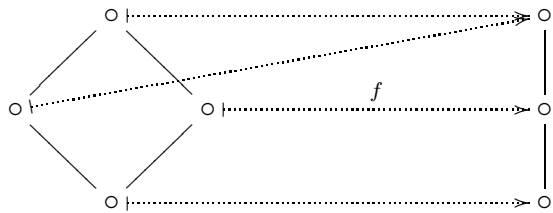
\[\xymatrix{&\circ\ar@{-}[dl]\ar@{-}[dr]\ar@{|.>}[rrrrr]&&&&&\circ\ar@{-}[d]\\
\circ\ar@{-}[dr]\ar@{|.>}[rrrrrru]&&\circ\ar@{-}[dl]\ar@{|.>}[rrrr]^f&&&&\circ\ar@{-}[d]\\
&\circ\ar@{|.>}[rrrrr]&&&&&\circ
}\]

\section{Adjointable maps}\label{sec:adjmaps}\label{sec:shhoms}

\subsection*{\sc Basic definitions and properties}

Similarly to Hilbert C*-modules, the module homomorphisms which have ``operator adjoints'' play a special role:

\begin{definition}
Let $X$ and $Y$ be pre-Hilbert $B$-modules. A function \[h:X\to Y\] is \emph{adjointable} if there is another function $h^\dagger:Y\to X$ such that for all $x\in X$ and $y\in Y$ we have
\[\inner{h(x)}y = \inner x{h^\dagger(y)}\;.\]
\end{definition}
[The usual notation for $h^\dagger$ is $h^*$, but we want to avoid confusion with the notation for inverse image homomorphisms of locale maps.]

Any adjointable function $h:X\to Y$ is necessarily a homomorphism of $B$-modules \cite{Paseka}, and this fact is a consequence of the non-degeneracy of $\inner--_Y$ alone; that is, $h$ satisfies $h\left(\V a_\alpha x_\alpha\right)=\V a_\alpha h(x_\alpha)$ because for all $y\in Y$ we have
\begin{eqnarray*}
\left\langle h\left(\V a_\alpha x_\alpha\right),y\right\rangle
&=&\left\langle \V a_\alpha x_\alpha,h^\dagger(y)\right\rangle
=\V a_\alpha\langle x_\alpha, h^\dagger(y)\rangle\\
&=&\V a_\alpha\langle h(x_\alpha),y\rangle
=\left\langle\V a_\alpha h(x_\alpha),y\right\rangle\;.
\end{eqnarray*}

\subsection*{\sc Hilbert bases}

The homomorphisms of Hilbert $B$-modules equipped with Hilbert bases are necessarily adjointable, and in order to prove this only the domain module need have a Hilbert basis:

\begin{theorem}
Let $X$ and $Y$ be pre-Hilbert $B$-modules such that $X$ has a Hilbert basis $\sections$ (hence, $X$ is Hilbert), and let $h:X\to Y$ be a homomorphism of $B$-modules. Then $h$ is adjointable with a unique adjoint $h^\dagger$, which is given by
\begin{eq}\label{adjoint}
h^\dagger(y) = \V_{t\in\sections}\langle h(t),y\rangle t\;.
\end{eq}
\end{theorem}

\begin{proof}
Let $x\in X$, $y\in Y$, and let us compute $\langle x,h^\dagger(y)\rangle$ using (\ref{adjoint}):
\begin{eqnarray*}
\langle x,h^\dagger(y)\rangle &=&\left\langle \V_{s\in\sections}\langle x,s\rangle s,\V_{t\in\sections}\langle h(t),y\rangle t\right\rangle\\
&=&
\V_{s,t\in\sections}\inner x s\wedge\inner s t\wedge\inner {h(t)}y\\
&=&
\V_{t\in\sections}\inner x t\wedge\inner {h(t)}y=
\left\langle \V_{t\in\sections}\inner x t h(t),y\right\rangle\\
&=&
\left\langle h\left(\V_{t\in\sections}\inner x t t\right),y\right\rangle=
\inner {h(x)}y\;.
\end{eqnarray*}
This shows that $h^\dagger$ is adjoint to $h$, and the uniqueness is a consequence of the non-degeneracy of the inner product of $X$. \qed
\end{proof}

\begin{corollary}\label{cor:homiffadj}
If $X$ and $Y$ are Hilbert $B$-modules and $X$ has a Hilbert basis then any function $h:X\to Y$ is adjointable if and only if it is a homomorphism of $B$-modules.
\end{corollary}

\begin{definition}
The \emph{category of Hilbert $B$-modules with Hilbert bases}, denoted by $B$-$\hmb$, is the category whose objects are those Hilbert $B$-modules for which there exist Hilbert bases and whose morphisms are the homomorphisms of $B$-modules (equivalently, the adjointable maps).
\end{definition}

\begin{corollary}
The assignment from homomorphisms $h$ to their adjoints $h^\dagger$ is a strong self-duality $(-)^\dagger:\opp{(B\textrm{-}\hmb)}\to B\textrm{-}\hmb$.
\end{corollary}

The matrix representations of Hilbert modules with Hilbert bases (\cf\ \ref{usefullemma}-\ref{ulmatrix} and \ref{lem:HMBisLOC}) can be extended to homomorphisms in a natural way. In particular the adjoint $h^\dagger$ of a homomorphism $h$ corresponds to the transpose of the matrix of $h$.
In order to see this consider the category of matrices $\Mat_B$ whose objects are the projection matrices over $B$ and whose arrows $F:M\to N$, given projection matrices
\begin{eqnarray*}
M&:&S\times S\to B\\
N&:&T\times T\to B\;,
\end{eqnarray*}
are the matrices $F:T\times S\to B$ such that $FM=F=NF$; the multiplication of arrows is just matrix multiplication, the unit arrows are the objects themselves, and this is a strongly self-dual category whose involution is transposition of matrices.

\begin{theorem}
The categories $B$-$\hmb$ and $\Mat_B$ are equivalent.
\end{theorem}

\begin{proof}
The assignment from modules $X$ to matrices
${\inner --}_X:\sections_X\times\sections_X\to B$
extends to the functor $\mathcal M:B\textrm{-}\hmb\to\Mat_B$ that to each homomorphism $h:Y\to X$ assigns the matrix $\mathcal M(h):\sections_X\times\sections_Y\to B$ defined by
\[(\mathcal M(h))_{st}=\inner{h(t)}s\;.\]
In the converse direction, the construction of a module $MB^S$ from each matrix $M:S\times S\to B$ extends to a functor $\mathcal X:\Mat_B\to B\textrm{-}\hmb$: given projection matrices
$M:S\times S\to B$ and $N:T\times T\to B$,
and an arrow $F:M\to N$, we define
$\mathcal X(F):M B^S\to NB^T$
by $\mathcal X(F)(f) = Ff$.
There is a natural isomorphism $\mathcal X\circ\mathcal M\cong \ident$, due to \ref{lem:HMBisLOC}, and a natural isomorphism $\mathcal M\circ\mathcal X\cong\ident$ follows from the equivalence of \cite{FS} between the category of $B$-sets and that of complete $B$-sets (see also \cite[\S2.9]{Borceux} or \cite[pp.\ 502--513]{Elephant}). Hence, the functors $\mathcal M$ and $\mathcal X$ form an adjoint equivalence of categories. \qed
\end{proof}

We remark that the maps of $B$-sets of \cite{FS} are arrows of $\Mat_B$, and thus the category of $B$-sets is a subcategory of $\Mat_B$.

\subsection*{\sc Sheaf homomorphisms}

In this section we exhibit an identification of ``operator adjoints'' with categorical adjoints, which as a consequence shows that the duality between homomorphisms of \'etale $B$-locales and sheaf homomorphisms is a restriction of the strong self-duality of $B$-$\hmb$.
In what follows we shall always consider an \'etale $B$-locale $X$ to be a Hilbert $B$-module with respect to the Hilbert basis of local sections $\sections_X$.

\begin{theorem}\label{lem:directimageisadjoint}
Let $X$ and $Y$ be \'etale $B$-locales, and let $f:X\to Y$ be a map of $B$-locales. Then $f_!=(f^*)^\dagger$ (equivalently, $f^*=(f_!)^\dagger$).
\end{theorem}

\begin{proof}
Let the $B$-locales, their projections, and $f$ be as follows:
\[
\xymatrix{
X\ar[rr]^f\ar[rd]_p&&Y\ar[ld]^q\\
&B
}
\]
Since $f$ commutes with the projections, which are local homeomorphisms, it is itself a local homeomorphism and thus it satisfies the Frobenius reciprocity condition $f_!(x\wedge f^*(y))=f_!(x)\wedge y$ (\ie, $f_!$ is $Y$-equivariant). Hence, we have, for all $x\in X$ and $y\in Y$:
\begin{eqnarray*}
\langle x,f^*(y)\rangle_X &=& p_!(x\wedge f^*(y))
=q_!(f_!(x\wedge f^*(y)))\\
&=&p_!(f_!(x)\wedge y)
=\langle f_!(x),y\rangle_Y\;. \qed
\end{eqnarray*}
\end{proof}

Now we shall look at a converse to the above theorem, whose proof depends on the following lemmas:

\begin{lemma}\label{lem:5.5}
Let $X$ be an \'etale $B$-locale. If $s\in \sections_X$ and $(b_\alpha)$ is a non-empty family of elements of $B$, we have $\left(\bigwedge_\alpha b_\alpha\right)s = \bigwedge_\alpha (b_\alpha s)$.
\end{lemma}

\begin{proof}
Let $(b_\alpha)$ be a non-empty family of elements of  $B$, and let $s\in\sections_X$. Then $\left(\bigwedge_\alpha b_\alpha\right)s$ is a lower bound of the set $\{b_\alpha \sigma\}$. Let $t$ be another lower bound. Then $\spp(t)\le b_\alpha$ for all $\alpha$ and, since $(b_\alpha)$ is non-empty, we have $t\le b_\alpha s$ for some $\alpha$ and thus $t\le s$. Hence,
\[\spp(t)\left(\left(\bigwedge_\alpha b_\alpha\right)s\right) =
    \left(\spp(t)\wedge\bigwedge_\alpha b_\alpha\right)s
    = \spp(t)s = t\;,\]
and it follows that $t\le\left(\bigwedge_\alpha b_\alpha\right)s$. This shows that $\left(\bigwedge_\alpha b_\alpha\right)s = \bigwedge_\alpha( b_\alpha s)$. \qed
\end{proof}

\begin{lemma}\label{lem:5.6}
Let $X$ be an \'etale $B$-locale, and let $S\subset\sections_X$ be a non-empty set such that $\V S\in\sections_X$. Then $\spp\left(\bigwedge S\right) = \bigwedge_{t\in S}\spp(t)$.
\end{lemma}

\begin{proof}
Let $s=\V S$. We have $\spp\left(\bigwedge S\right)
= \spp\left(\bigwedge_{t\in S}\spp(t)s\right)$ and, by \ref{lem:5.5}, this equals
\[\spp\left(\left(\bigwedge_{t\in S}\spp(t)\right)s\right)=\left(\bigwedge_{t\in S}\spp(t)\right)\wedge\spp\left(s\right)\;.\]
Since $S$ is non-empty the latter equals $\bigwedge_{t\in S}\spp(t)$. \qed
\end{proof}

\begin{theorem}\label{lem:adjointsandmeets}
Let $h:X\to Y$ be a sheaf homomorphism of \'etale $B$-locales. Then its adjoint $h^\dagger$ preserves arbitrary meets.
\end{theorem}

\begin{proof}
Let $S\subset Y$. We shall show that $h^\dagger\left(\bigwedge S\right)=\bigwedge h^\dagger(S)$ by using the non-degeneracy of the inner product of $X$; that is, we shall prove, for all $s\in\sections_X$, that
$\left\langle s,h^\dagger\left(\bigwedge S\right)\right\rangle=\left\langle s,\bigwedge h^\dagger(S)\right\rangle$.
Let then $s\in\sections_X$. We have
\begin{eqnarray*}
\left\langle s,h^\dagger\left(\bigwedge S\right)\right\rangle &=&
\left\langle h(s),\bigwedge S\right\rangle = \spp\left(h(s)\wedge\bigwedge S\right)\\
&=& \spp\left(h(s)\wedge\bigwedge_{y\in S} (h(s)\wedge y)\right) =\spp\left(\bigwedge S'\right)\;,
\end{eqnarray*}
where the set $S'=\{h(s)\}\cup\{h(s)\wedge y\st y\in S\}$ is non-empty,
it is contained in $\sections_Y$ because $h(s)\in\sections_Y$, and $\V S'\in\sections_Y$ because $S'$ is upper bounded by $h(s)$. Hence, by \ref{lem:5.6}, we have $\spp\left(\bigwedge S'\right)=\bigwedge\spp(S')$. Moreover, $\spp(h(s))=\spp(s)$ and thus
\begin{eqnarray*}
\bigwedge\spp(S') &=& \spp(h(s))\wedge\bigwedge_{y\in S}\spp(h(s)\wedge y)
=\spp(s)\wedge\bigwedge_{y\in S}\langle h(s),y\rangle\\
&=&\spp(s)\wedge\bigwedge_{y\in S}\langle s,h^\dagger(y)\rangle
=\spp(s)\wedge\bigwedge_{y\in S}\spp(s\wedge h^\dagger(y))
=\bigwedge\spp(S'')\;,
\end{eqnarray*}
where the set $S''=\{s\}\cup\{s\wedge h^\dagger(y)\st y\in S\}$ is non-empty, it is contained in $\sections_X$, and $\V S''\in\sections_X$ because $S''$ is upper bounded by $s$. Hence, again by \ref{lem:5.6}, we have
\begin{eqnarray*}
\bigwedge\spp(S'') &=& \spp\left(\bigwedge S''\right)
= \spp\left(s\wedge\bigwedge_{y\in S} (s\wedge h^\dagger(y))\right)\\
&=& \spp\left(s\wedge\bigwedge h^\dagger(S)\right)
=\left\langle s,\bigwedge h^\dagger(S)\right\rangle\;,
\end{eqnarray*}
which concludes the proof. \qed
\end{proof}

Theorem \ref{thm:postponed}, whose proof has been postponed until now, is a simple corollary of these results. Instead of having proved directly that the right adjoints of sheaf homomorphisms are module homomorphisms, as we might have attempted in \S\ref{sec:shhoms1}, we have instead shown that the ``operator adjoints'' of sheaf homomorphisms, which are module homomorphisms, are also homomorphisms of locales:

\begin{corollary}[\cf\ Theorem \ref{thm:postponed}]
The categories $B$-$\sh$ and $B$-$\Etale$ are isomorphic.
\end{corollary}

\begin{proof}
By \ref{lem:adjointsandmeets}, the adjoint $h^\dagger$ of a sheaf homomorphism $h$ is a homomorphism of $B$-locales. This defines a map of $B$-locales $f$ such that $f^*=h^\dagger$. By \ref{lem:directimageisadjoint}, $f_!=(f^*)^\dagger=h$, and thus the faithful functor $\mathcal S$ of \ref{thm:postponed} is full. \qed
\end{proof}

~\\
{\sc
Departamento de Matem\'{a}tica\\
Instituto Superior T\'{e}cnico\\
Universidade T\'{e}cnica de Lisboa\\
Av.\ Rovisco Pais 1, 1049-001 Lisboa, Portugal\\
{\it E-mail:} {\sf pmr@math.ist.utl.pt}\\
~\\
Departamento de Matem\'{a}tica e Engenharias\\
Campus Universit\'{a}rio da Penteada\\
Universidade da Madeira\\
9000-390 Funchal, Portugal\\
{\it E-mail:} {\sf elias@uma.pt}
}

\begin{thebibliography}{10}
{\small
\bibitem{Borceux} F.\ Borceux, Handbook of Categorical Algebra 3 --- Categories of Sheaves, Encyclopedia of Math.\ Appl., vol.\ 52, Cambridge Univ.\ Press, 1994.

\bibitem{BB} F.\ Borceux, G.\ van den Bossche, Quantales and their sheaves, Order 3 (1986) 61--87.

\bibitem{FS} M.P.\ Fourman, D.S.\ Scott, Sheaves and logic, In: M.P.\ Fourman, C.J.\ Mulvey, D.S.\ Scott (Eds.), Applications of Sheaves, Lect.\ Notes in Math., Vol.\ 753, Springer, 1979, pp.\ 302--401.

\bibitem{G00}
R.P.\ Gylys, Sheaves on quantaloids, Lithuanian Mathematical Journal 40 (2000) 105--134.

\bibitem{G01}
R.P.\ Gylys, Sheaves on involutive quantaloids, Lithuanian Mathematical Journal 41 (2001) 35--53.

\bibitem{KR} D.\ Kruml, P.\ Resende, On quantales that classify C*-algebras, Cah.\ Topol.\ G\'eom.\ Diff\'er.\ Cat\'eg.\ 45 (2004) 287--296.

\bibitem{stone}
P.T.\ Johnstone,
Stone Spaces,
Cambridge Stud.\ Adv.\ Math.,
vol.\ 3, Cambridge Univ.\ Press, 1982.

\bibitem{Elephant} P.T.\ Johnstone, Sketches of an Elephant: A Topos Theory Compendium (Vol.\ 2), Oxford Logic Guides, No.\ 44, Clarendon Press, 2002.

\bibitem{JT} A.\ Joyal, M.\ Tierney, An Extension of the Galois Theory of Grothendieck, Mem.\ Amer.\ Math.\ Soc., Vol.\ 309, American Mathematical Society, 1984.

\bibitem{Lance} E.C.\ Lance, Hilbert C*-Modules---A Toolkit for Operator Algebraists, London Mathematical Society Lecture Note Series, No.\ 210, Cambridge University Press, 1995.

\bibitem{Mulvey}
C.J.\ Mulvey,
Quantales,
In: M.\ Hazewinkel (Ed.), The Encyclopaedia of Mathematics, third supplement, Kluwer Academic Publishers, 2002, pp.\ 312--314.

\bibitem{MN}
C.J.\ Mulvey, M.\ Nawaz, Quantales: Quantal sets, In: Non-Classical Logics and their Application to Fuzzy Subsets: A Handbook of the Mathematical Foundations of Fuzzy Set Theory, Kluwer, 1995, pp.\ 159--217.

\bibitem{Paseka} J.\ Paseka, Hilbert $Q$-modules and nuclear ideals in the category of $\bigvee$-semilattices with a duality. CTCS '99: Conference on Category Theory and Computer Science (Edinburgh), Paper No.\ 29019, 19 pp.\ (electronic), Electron.\ Notes Theor.\ Comput.\ Sci., 29, Elsevier, Amsterdam, 1999.

\bibitem{Paseka2} J.\ Paseka, Hermitian kernels, Hilbert $Q$-modules and Ando dilation. Contributions to general algebra, 12 (Vienna, 1999), 317--335, Heyn, Klagenfurt, 2000.

\bibitem{PR}
J.\ Paseka, J.\ Rosick\'{y},
Quantales,
In: B.\ Coecke, D.\ Moore, A.\ Wilce, (Eds.), Current
Research in Operational Quantum Logic: Algebras,
  Categories and Languages,
Fund.\ Theories
  Phys., vol.\ 111, Kluwer Academic Publishers, 2000,
pp.\ 245--262.

\bibitem{JoelPhD}
J.Z.\ Ramos, Involutive Quantales, PhD Thesis, Univ.\ Sussex, 2006.

\bibitem{AIM} P.\ Resende, \'Etale groupoids and their quantales, Adv.\ Math.\ 208 (2007) 147--209.

\bibitem{am} P.\ Resende, Quantal sets, quantale modules, and groupoid actions, talk at Internat.\ Conf.\ Category Theory 2007, June 17--23, 2007, Carvoeiro, Portugal; http://www.mat.uc.pt/$\sim$categ/ct2007/slides/resende.pdf (paper in preparation).

\bibitem{RV} P.\ Resende, S.J.\ Vickers, Localic sup-lattices and tropological systems, Theoret.\ Comput.\ Sci.\ 305 (2003) 311--346.

\bibitem{Rosenthal}
K.\ Rosenthal,
Quantales and Their Applications, Pitman Research Notes in Mathematics Series 234, Longman Scientific \& Technical, 1990.

\bibitem{Rosenthal2}
K.\ Rosenthal,
The Theory of Quantaloids, Pitman Research Notes in Mathematics Series 348, Addison Wesley Longman Limited, 1996.

\bibitem{IsarPhD}
I.\ Stubbe, Categorical structures enriched in a quantaloid: categories and semicategories, PhD Thesis, Univ.\ Louvain-la-Neuve, 2003.

\bibitem{S07}
I.\ Stubbe, $Q$-modules are $Q$-suplattices, Theory and Applications of Categories 19 (2007) 50--60.

\bibitem{W}
R.F.C.\ Walters, Sheaves on sites as Cauchy complete categories, 
J.\ Pure and Applied Algebra 24 (1982) 95--102.
}
\end{thebibliography}
\end{document}